\documentclass[12pt]{article}
\usepackage{epsfig, subfigure, verbatim}
\usepackage{amssymb}

\author{
Vladimir E. Alekseev\thanks{Department of Mathematical Logic, University of Nizhny Novgorod, Gagarina 23, 603950 RUSSIA. E-mail: ave@uic.nnov.ru},\\ Alastair Farrugia\thanks{Department of Combinatorics and Optimization, University of Waterloo, Waterloo, Ontario, Canada, N2L6G1. Email: afarrugia@math.uwaterloo.ca}, \\ 
Vadim V. Lozin\thanks{RUTCOR, Rutgers University, 640 Bartholomew Road, Piscataway, NJ 08854-8003, USA. E-mail: lozin@rutcor.rutgers.edu}}
\title{New Results on Generalized Graph Coloring}
\date{}
\newcommand{\cK}{{\cal K}}
\newcommand{\cO}{{\cal O}}
\newcommand{\cP}{{\cal P}}
\newcommand{\cQ}{{\cal Q}}

\newtheorem{thm}{Theorem}

\newtheorem{lem}{Lemma}
\newenvironment{proof}[1][Proof]{\noindent\textbf{#1.} }{\ \rule{0.5em}{0.5em}}

\begin{document}
\maketitle

\medskip
\begin{abstract} For graph classes $\cP_1,\ldots,\cP_k$, Generalized Graph Coloring is the problem of deciding whether the vertex set of a given graph $G$ can be partitioned into subsets $V_1,\ldots,V_k$ so that $V_j$ induces a graph in the class $\cP_j$ $(j=1,2,\ldots,k)$. If $\cP_1=\cdots=\cP_k$ is the class of edgeless graphs, then this problem coincides with the standard vertex $k$-{\sc colorability}, which is known to be NP-complete for any $k\ge 3$. Recently, this result has been generalized by showing that if all $\cP_i$'s are 
%CHANGED
additive induced-hereditary, 
then the generalized graph coloring is NP-hard, with the only exception of bipartite graphs. Clearly, a similar result follows when all the $\cP_i$'s are co-additive. 

In this paper, we study the problem where we have a mixture of additive and co-additive classes, presenting several new results dealing both with NP-hard and polynomial-time solvable instances of the problem. 
%+%
%In fact, our algorithm works for any fixed graphs $G_1, \ldots, G_r, H_1, %\ldots, H_s$, and positive integers $n$ and $m$, recognising $(Free(K_n, G_1, %\ldots, G_r), Free(\overline{K_m}, H_1, \ldots, H_s)$-colorable graphs in %polynomial-time.
\end{abstract}


%%%%%%%%%%%%%%%%%%%%%%%%%%%%%%%%%%%%%%%%%%%%%%%%%%%%%%%%%%%%%%%%%%%%%
%%%%%%%%%%%%%%%%%%%%%%%%%%%%%%%%%%%%%%%%%%%%%%%%%%%%%%%%%%%%%%%%%%%%%
%%%%%%%%%%%%%%%%%%%%%%%%%%%%%%%%%%%%%%%%%%%%%%%%%%%%%%%%%%%%%%%%%%%%%
%%%%%%%%%%%%%%%%%%%%%%%%%%%%%%%%%%%%%%%%%%%%%%%%%%%%%%%%%%%%%%%%%%%%%

\section{Introduction}
All graphs in this paper are finite, without loops and multiple edges. For a graph $G$ we denote by $V(G)$ and $E(G)$ the vertex set and the edge set of $G$, respectively. By $N(v)$ we denote the neighborhood of a vertex $v\in V(G)$, i.e. the subset of vertices of $G$ adjacent to $v$. The subgraph of $G$ induced by a set $U\subseteq V(G)$ will be denoted $G[U]$. We say that a graph $G$ is $H$-free if $G$ does not contain $H$ as an induced subgraph. As usual, $K_n$ and $P_n$ stand for the complete graph and chordless path on $n$ vertices, respectively, and the complement of a graph $G$ is denoted $\overline{G}$.

A class of graphs, or synonymously graph property, $\cP$ is said to be {\it hereditary} if $G\in \cP$ implies $G-v\in \cP$ for any vertex $v\in V(G)$. We call $\cP$ {\it monotone} if $G\in \cP$ implies $G-v\in \cP$ for any vertex $v\in V(G)$ and $G-e\in \cP$ for any edge $e\in E(G)$. Clearly every monotone property is hereditary, but the converse statement is not true in general. A property $\cP$ is {\it additive} if $G_1\in \cP$ and $G_2\in \cP$ with $V(G_1)\cap V(G_2)=\emptyset$ implies $G=(V(G_1)\cup V(G_2),E(G_1)\cup E(G_2))\in \cP$. The class of graphs containing no induced subgraphs isomorphic to graphs in a set $Y$ will be denoted $Free(Y)$.  It is well known that a class of graphs $\cP$ is hereditary if and only if $\cP=Free(Y)$ for some set $Y$. 

A property is said to be non-trivial if it contains at least one, but not all graphs. The {\em complementary property\/} of $\cP$ is 
$\overline{\cP} := \{\overline{G} \mid G \in \cP\}$.
Note that $\cP$ is he\-re\-di\-ta\-ry if and only if $\overline{\cP}$ is.
So a {\em co-additive he\-re\-di\-ta\-ry\/} property, i.e. the complement of an additive he\-re\-di\-ta\-ry property, is itself he\-re\-di\-ta\-ry.

Let $\cP_1,\ldots,\cP_k$ be graph properties (classes) with $k>1$. A graph $G=(V,E)$ is  {\em $(\cP_1,$ $\ldots,$ $\cP_k)$-colorable\/} if there is a partition $(V_1,\ldots,V_k)$ of $V(G)$ such that $G[V_j]\in \cP_j$ for each $j=1,\ldots,k$. The problem of recognizing $(\cP_1,$ $\ldots,$ $\cP_k)$-colorable graphs is usually referred to as Generalized Graph Coloring \cite{Bro96}. When $\cP_1=\cdots=\cP_k$ is the class $\cO$ of edgeless graphs, this problem coincides with the standard $k$-{\sc colorability}, which is known to be NP-complete for $k\ge 3$. Generalized Graph Coloring remains difficult for many other cases. For example, Cai and Corneil~\cite{cai} showed that ($Free(K_n)$,$Free(K_m)$)-coloring is NP-complete for any integers $m,n \geq 2$, with the exception $m=n=2$. 
%CHANGED These results as well as some other results on this topic 
This result, and others
\cite{Ach97,Bro96,KS97}, have been recently generalized in \cite{farr:comp} as follows.

\begin{thm}\label{thm:add}
If $\cP_1,\ldots,\cP_k$ ($k>1$) are additive hereditary classes of graphs, then the problem of recognizing $(\cP_1,$ $\ldots,$ $\cP_k)$-colorable graphs is NP-hard, unless $k=2$ and $\cP_1=\cP_2$ is the class of edgeless graphs.  
\end{thm}

Clearly, a similar result follows for co-additive properties. In the present paper we focus on the case where we have a mixture of additive and co-additive properties. 
%+%
%In fact, our algorithm works for any fixed graphs $G_1, \ldots, G_r, H_1, %\ldots, H_s$, and positive integers $n$ and $m$, recognising $(Free(K_n, G_1, %\ldots, G_r), Free(\overline{K_m}, H_1, \ldots, H_s)$-colorable graphs in %polynomial-time.

The {\em product\/} of graph classes $\cP_1, \ldots, \cP_k$ is $\cP_1 \circ \cdots \circ \cP_k := \{G \mid G$ is $(\cP_1, \ldots, \cP_k)$-colorable$\}$. A property is {\em reducible\/} if it is the product of two other properties, otherwise it is {\em irreducible}. It can be easily checked that the product of additive hereditary (or monotone) properties is again additive hereditary (respectively, monotone); and that $\overline{\cP_1 \circ \cdots \circ \cP_k } = \overline{\cP_1} \circ  \cdots \circ \overline{\cP_k}$. So, without loss of generality we shall restrict our study to the case $k=2$ and shall denote throughout the paper an additive property by $\cP$ and co-additive by $\cQ$. We will refer to the problem of recognizing $(\cP, \cQ)$-colorable graphs as $(\cP \circ \cQ)$-\textsc{recognition}. 

The plan of the paper is as follows. In Section~\ref{sec:hard}, we show that $(\cP \circ \cQ)$-\textsc{recognition} cannot be simpler than $\cP$- or  $\cQ$-\textsc{recognition}. 
%CHANGED Specifically, 
In particular,
we prove that $(\cP \circ \cQ)$-\textsc{recognition} is NP-hard whenever $\cP$- or  $\cQ$-\textsc{recognition} is NP-hard. Then, in Section~\ref{sec:pol}, we study the problem under the assumption that both $\cP$- and $\cQ$-\textsc{recognition} are polynomial-time solvable and present infinitely many classes of $(\cP, \cQ)$-colorable graphs with polynomial recognition time. These two results together give a complete answer to the question of complexity of $(\cP \circ \cQ)$-\textsc{recognition} when $\cP$ and $\overline{\cQ}$ are additive monotone. When $\cP$ and $\overline{\cQ}$ are additive hereditary (but not both monotone), there remains an unexplored gap that we discuss in the concluding section of the paper. 

%%%%%%%%%%%%%%%%%%%%%%%%%%%%%%%%%%%%%%%%%%%%%%%%%%%%
%%%%%%%%%%%%%%%%%%%%%%%%%%%%%%%%%%%%%%%%%%%%%%%%%%%%
\section{NP-hardness}
\label{sec:hard}
%%%%%%%%%%%%%%%%%%%%%%%%%%%%%%%%%%%%%%%%%%%%%%%%%%%%
%%%%%%%%%%%%%%%%%%%%%%%%%%%%%%%%%%%%%%%%%%%%%%%%%%%%

In this section we prove that if $\cP$-\textsc{recognition} (or 
%COMMENTED OUT equivalently 
$\cQ$-\textsc{recognition}) is NP-hard, then so is 
$(\cP \circ \cQ)$-\textsc{recognition}. This is a direct consequence of the theorem below. In this theorem we use uniquely colorable graphs, which are often a crucial tool in proving coloring results.

A graph $G$ is {\em uniquely $(\cP_1, \ldots, \cP_k)$-colorable\/} if $(V_1, \ldots, V_k)$ is its only $(\cP_1, \ldots, \cP_k)$-partition, up to some permutation of the $V_i$'s. 
If, say, $\cP_1 = \cP_2$, then $(V_2, V_1, V_3$,$\ldots$,$V_k)$ will also be a $(\cP_1, \cP_2, \cP_3, \ldots, \cP_k)$-coloring of $G$; such a permutation (of $V_i$'s that correspond to equal properties) is a {\em trivial interchange\/}.
A graph is {\em strongly uniquely $(\cP_1, \ldots, \cP_k)$-colorable\/} if 
$(V_1, \ldots, V_k)$ is the only $(\cP_1, \ldots, \cP_k)$-coloring, up to trivial interchanges. 

When $\cP_1, \ldots, \cP_k$ are 
%CHANGED irreducible additive he\-re\-di\-ta\-ry properties, 
irreducible he\-re\-di\-ta\-ry properties, and each $\cP_i$ is either additive or co-additive,
there is a strongly uniquely $(\cP_1, \ldots, \cP_k)$-colorable graph with each $V_i$ non-empty. This important construction, for additive $\cP_i$'s, is due to Mih\'ok~\cite{uft-2}, with some embellishments by Broere and Bucko~\cite{uni-1}, while the proof of unique colorability follows from~\cite[Thm. 5.3]{discussiones}. Obviously, similar results apply to co-additive properties. The generalization to mixtures of additive and co-additive properties can be found in~\cite[Cor. 4.3.6, Thm. 5.3.2]{thesis} 

\begin{thm}\label{thm:hard}
Let $\cP$ and $\overline{\cQ}$ be additive he\-re\-di\-ta\-ry properties. Then there is a poly\-no\-mial-time reduction from $\cP$-\textsc{recognition} to 
$(\cP \circ \cQ)$-\textsc{recognition}. 
\end{thm}

\begin{proof}
Let $\cP = \cP_1 \circ \cdots \circ \cP_n$ and $\cQ = \cQ_{1} \circ \cdots \circ \cQ_{r}$, where the $\cP_i$'s and $\overline{\cQ}_j$'s are irreducible additive hereditary properties. As noted above, there is a strongly uniquely
$(\cP_1, \ldots, \cP_{n}, \cQ_1, \ldots, \cQ_r)$-colorable graph $H$ with
partition $(U_1, \ldots, U_{n}, W_1,$ $\ldots,$ $W_r)$, where each $U_i$ and $W_j$ is non-empty. Define $U := U_1 \cup \cdots \cup U_n$ and $W := W_1 \cup \cdots \cup W_r$. Arbitrarily fix a vertex $u \in U_1$, and define $N_W(u) := N(u) \cap W$. For any graph $G$, let the graph $G_H$ consist of disjoint copies of $G$ and $H$, together with edges $\{vw \mid v \in V(G), w \in N_W(u)\}$.
We claim that $G_H \in \cP \circ \cQ$ if and only if $G \in \cP$.

\begin{figure}[htb]
\begin{center}
\input{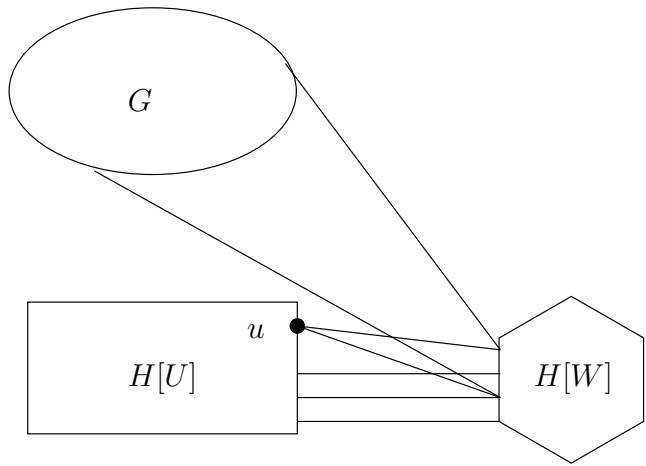}
\caption{Using $H$ to construct $G_H$.}
\end{center}
\end{figure}

If $G \in \cP$, then, by additivity, $G \cup H[U]$ is in $\cP$, and thus $G_H$ is in $\cP \circ \cQ$. Conversely, suppose $G_H \in \cP \circ \cQ$, i.e. it has a $(\cP_1, \ldots, \cP_n, \cQ_{1}, \ldots, \cQ_{r})$-partition, say $(X_1, \ldots, X_n, Y_1, \ldots, Y_r)$. Since $H$ is strongly uniquely partitionable, we can assume that, for $1 \leq i \leq r$, $Y_i \cap V(H) = W_i$. Now, suppose for contradiction that, for some $k$, there is a vertex $v \in V(G)$ such that $v \in Y_k$; without loss of generality, let $k = r$. Then $G_H[W_r \cup \{v\}] \cong H[W_r \cup \{u\}]$ is in $\cQ_r$, so $(U_1 \setminus \{u\}, U_2, \ldots, U_n, W_1, \ldots, W_{r-1}, W_r \cup \{u\})$ is a new 
$(\cP_1, \ldots, \cP_{n}, \cQ_1, \ldots, \cQ_r)$-partition of $H$, which is impossible.
Thus, $V(G) \subseteq X_1 \cup \cdots \cup X_n$, and hence $G \in \cP$, as claimed.

Since $H$ is a fixed graph, $G_H$ can be constructed in time linear in $|V(G)|$, so the theorem is proved.
\end{proof}

%%%%%%%%%%%%%%%%%%%%%%%%%%%%%%%%%%%%%%%%%%%%%
%%%%%%%%%%%%%%%%%%%%%%%%%%%%%%%%%%%%%%%%%%%%%
\section{Polynomial time results}
\label{sec:pol}
%%%%%%%%%%%%%%%%%%%%%%%%%%%%%%%%%%%%%%%%%%%%%
%%%%%%%%%%%%%%%%%%%%%%%%%%%%%%%%%%%%%%%%%%%%%

\begin{lem}\label{lem:pol}
For any $\cP\subseteq Free(K_n)$ and $\cQ\subseteq Free(\overline{K}_m)$, there exists a constant $\tau=\tau(\cP,\cQ)$ such that for every graph $G=(V,E)\in \cP\circ\cQ$ and every subset $B\subseteq V$ with $G[B]\in \cP$, at least one of the following statements holds:
\begin{itemize}
\item[(a)] there is a subset $A\subseteq V$ such that $G[A]\in \cP$, $G[V-A]\in \cQ$, and $|A-B|\le \tau$,
\item[(b)] there is a subset $C\subseteq V$ such that $G[C]\in \cP$, $|C|=|B|+1$, and
$|B-C|\le \tau$.
\end{itemize}
\end{lem}

\begin{proof}
By the Ramsey Theorem \cite{GRS80}, for each positive integers $m$ and $n$, there is a constant $\tau(m,n)$ such that every graph with more than $\tau(m,n)$ vertices contains either a $\overline{K}_m$ or a $K_n$ as an induced subgraph. For two classes $\cP\subseteq Free(K_n)$ and $\cQ\subseteq Free(\overline{K}_m)$, we define $\tau=\tau(\cP,\cQ)$ to be equal $\tau(m,n)$. Let us show that with this definition the proposition follows.

Let $G=(V,E)$ be a graph in $\cP\circ\cQ$, and $B$ a subset of $V$ such that $G[B]\in \cP$. Consider an arbitrary subset $A\subseteq V$ such that $G[A]\in \cP$ and $G[V-A]\in \cQ$. If (a) does not hold, then $|A-B|>\tau$. Furthermore, $G[B-A]\in \cP\cap \cQ\subseteq Free(K_n,\overline{K}_m)$, and hence $|B-A|\le \tau$. Therefore, $|A|>|B|$. But then any subset $C\subseteq A$ such that $A\cap B\subseteq C$ and $|C|=|B|+1$ satisfies (b).  
\end{proof}

\medskip
Lemma~\ref{lem:pol} suggests the following recognition algorithm for graphs in the class $\cP\circ\cQ$.  

\medskip
{\bf Algorithm $\cal{A}$}

\medskip
{\bf Input:} A graph $G=(V,E)$.

{\bf Output:} {\bf YES} if $G\in \cP\circ\cQ$, or {\bf NO} otherwise.
\begin{itemize}
\item[(1)]  Find in $G$ any 
%CHANGED maximal under inclusion 
inclusion-wise maximal
subset $B\subseteq V$ inducing a $K_n$-free graph.
\item[(2)] If there is a subset $C\subseteq V$ satisfying condition (b) of Lemma~\ref{lem:pol},\\ then set $B:=C$ and repeat Step (2). 
\item[(3)] If $G$ contains a subset $A\subseteq V$ such that
\begin{itemize}
\item[] $|B-A|\le \tau$,
\item[] $|A-B|\le \tau$,
\item[] $G[A]\in \cP$,
\item[] $G[V-A]\in \cQ$,
\end{itemize}
output {\bf YES}, otherwise output {\bf NO}.
\end{itemize}

\begin{thm}\label{thm:pol}
If graphs on $p$ vertices in a class $\cP\subseteq Free(K_n)$ can be recognized in time $O(p^k)$ and graphs in a class $\cQ\subseteq Free(\overline{K}_m)$ can be recognized in time $O(p^l)$, then Algorithm $\cal{A}$ recognizes graphs on $p$ vertices in the class $\cP \circ \cQ$ in time $O(p^{2\tau +\max\{(k+2),\max\{k,l\}\}})$, where $\tau=\tau(\cP,\cQ)$.
\end{thm}

\begin{proof}
Correctness of the algorithm follows from Lemma~\ref{lem:pol}. Now let us estimate its time complexity. In Step (2), the algorithm examines at most $p\choose\tau$$p\choose\tau+1$ subsets $C$ and for each of them verifies whether $G[C]\in \cP$ in time $O(p^k)$. Since Step (2) loops at most $p$ times, its time complexity is $O(p^{2\tau+k+2})$. In Step (3), the algorithm examines at most $p\choose\tau$$^2$ subsets $A$, and for each $A$, it verifies whether $G[A]\in \cP$ in time $O(p^k)$ and whether $G[V-A]\in \cQ$ in time $O(p^{l})$. Summarizing, we conclude that the total time complexity of the algorithm is $O(p^{2\tau +\max\{(k+2),\max\{k,l\}\}})$.
\end{proof}

\medskip
Notice that Theorem~\ref{thm:pol} generalizes several positive results on the topic under consideration. For instance, the split graphs \cite{FH77}, which are ($Free(K_2)$,$Free(\overline{K_2})$)-colorable by definition, can be recognized in polynomial time. More general classes have been studied under the name of polar graphs in \cite{CC86,MK85,TC85}. By definition, a graph is $(m-1,n-1)$ polar if it is $(\cP,\cQ)$-colorable with $\cP=Free(K_n,P_3)$ and $\cQ=Free(\overline{K}_m,\overline{P}_3)$. 
%+%
It is shown in \cite{MK85} that 
for any particular values of $m\ge 2$ and $n\ge 2$, $(m-1,n-1)$ polar graphs on $p$ vertices can be recognized in time $O(p^{2m+2n+3})$. 

Further examples generalizing the split graphs were examined in \cite{BLS98}
and~\cite{feder}, where the authors showed that classes of graphs partitionable into at most two independent sets and two cliques can be recognized in polynomial time. These are special cases of  $(\cP \circ \cQ)$-\textsc{recognition} with $\cP\subseteq Free(K_3)$ and $\cQ\subseteq Free(\overline{K}_3)$. 

%%%%%%%%%%%%%%%%%%%%%%%%%%%%%%%%%%%%%%%%%%%%%
%%%%%%%%%%%%%%%%%%%%%%%%%%%%%%%%%%%%%%%%%%%%%
\section{Concluding results and open problems} 
%%%%%%%%%%%%%%%%%%%%%%%%%%%%%%%%%%%%%%%%%%%%%
%%%%%%%%%%%%%%%%%%%%%%%%%%%%%%%%%%%%%%%%%%%%%
Theorems~\ref{thm:hard} and ~\ref{thm:pol} together provide complete answer to the question of complexity of $(\cP \circ \cQ)$-\textsc{recognition} in case of monotone properties $\cP$ and $\overline{\cQ}$. Indeed, if $\cP$ is an additive monotone non-trivial property, then $\cP\subseteq Free(K_n)$ for a certain value of $n$, since otherwise it includes all graphs. Similarly, if $\overline{\cQ}$ is additive monotone, then $\cQ\subseteq Free(\overline{K}_m)$ for some $m$. Hence, the following theorem holds. 

\begin{thm}\label{co-her-additive}
If $\cP$ and $\overline{\cQ}$ are additive monotone properties, then 
$(\cP \circ {\cQ})$-\textsc{re\-co\-gni\-tion} has poly\-no\-mial-time complexity if and only if $\cP$- and $\cQ$-\textsc{re\-co\-gni\-tion} are both poly\-no\-mial-time solvable; moreover, $(\cP \circ {\cQ})$-\textsc{recognition} is in NP if and only if $\cP$- and $\cQ$-\textsc{recognition} are both in NP. 
\end{thm}
 
If $\cP$ and $\overline{\cQ}$ are general additive hereditary properties (not necessarily monotone), then there is an unexplored gap containing properties $\cP \circ \cQ$, where $\cP$ and $\cQ$ can both be recognized in polynomial time, but $\cK\subset \cP$ or $\cO\subset \cQ$ 
%ADDED
(where $\cK := \overline{\cO}$ is the set of cliques). 
In the rest of this section we show that this gap contains both NP-hard and polynomial-time solvable instances, and propose several open problems to study. 

For a polynomial time result we refer the reader to \cite{TC85}, where the authors claim that $(\cP \circ \cQ)$-\textsc{recognition} is polynomial-time solvable if $\cP$ is the class of edgeless graphs and $\cQ=Free(\overline{P}_3)$. Notice that $Free(\overline{P}_3)$ contains all edgeless graphs and hence Theorem~\ref{thm:pol} does not apply to this case. Interestingly enough, when we extend $\cP$ to the class of bipartite graphs, we obtain an NP-hard instance of the problem, as the following theorem shows.

\begin{thm}\label{thm:ex1}
If $\cP$ is the class of bipartite graphs and $\cQ=Free(\overline{P}_3)$, then $(\cP \circ \cQ)$-\textsc{recognition} is NP-hard.
\end{thm}

\begin{proof}
We reduce the standard 3-{\sc colorability} to our problem. Consider an arbitrary graph $G$ and let $G'$ be the graph obtained from $G$ by adding a triangle $T=(1,2,3)$ with no edges between $G$ and $T$. We claim that $G$ is 3-colorable if and only if $G'$ is $(\cP, {\cQ})$-{colorable}.

First, assume that $G$ is 3-colorable and let $V_1,V_2,V_3$ be a partition of $V(G)$ into three independent sets. We define $V'_j=V_j\cup \{j\}$ for $j=1,2,3$. Then $G'[V'_1\cup V'_2]$ is a bipartite graph and $G'[V'_3]\in Free(\overline{P}_3)$, and the proposition follows.

Conversely, let $U\cup W$ be a partition of $V(G')$ with $G'[U]$ being a bipartite graph and $G'[W]\in Free(\overline{P}_3)$. Clearly, $T\nsubseteq U$. 
If $T-U$ contains a single vertex, then $G'[W-T]$ is an edgeless graph, since otherwise a $\overline{P}_3$ arises. If $T-U$ contains more than one vertex, then $W-T=\emptyset$ for the same reason. Clearly, in both cases
$G$ is a 3-colorable graph. 
\end{proof}

\medskip
This discussion presents the natural question of exploring the boundary that separates polynomial from non-polynomial time solvable instances in the above-mentioned gap. As one of the smallest classes in this gap with unknown recognition time complexity, let us point out $(\cP, {\cQ})$-\textsc{colorable} graphs with $\cP=\cO$ and $\cQ=Free(2K_2,P_4)$, where $2K_2$ is the disjoint union of two copies of $K_2$. 

Another direction for prospective research deals with $(\cP,\cQ)$-colorable graphs where $\cP$ or $\cQ$ is neither additive nor co-additive. This area seems to be almost unexplored and also contains both NP-hard and polynomial-time solvable problems. To provide some examples, let $\cQ$ be the class of complete bipartite graphs, which is obviously neither additive nor co-additive. The class of graphs partitionable into an independent set and a complete bipartite graph has been studied in \cite{BHLL02} under the name of bisplit graphs and has been shown there to be polynomial-time recognizable. Again, extension of $\cP$ to the class of all bipartite graphs transforms the problem into an NP-hard instance.    

\begin{thm}\label{thm:ex2}
If $\cP$ is the class of bipartite graphs and $\cQ$ is the class of complete bipartite graphs, then $(\cP \circ \cQ)$-\textsc{recognition} is NP-hard.
\end{thm}

\begin{proof}
The reduction is again from 3-{\sc colorability}. For a graph $G$, we define $G'$ to be the graph obtained from $G$ by adding a new vertex adjacent to every vertex of $G$. It is a trivial exercise to verify that $G$ is 3-colorable if and only if $G'$ is $(\cP, {\cQ})$-\textsc{colorable}.  
\end{proof}

\section{Acknowledgements}
The second author wishes to thank R.\ Bruce Richter, his doctoral supervisor, for his valuable comments; as well as the Canadian government, which is fully funding his studies in Waterloo through a Canadian Commonwealth Scholarship.

\end{document}